\documentclass[11pt,reqno]{amsart}
\usepackage{amsmath,amssymb,graphicx}
\newtheorem{thm}{Theorem}[section]
\newtheorem{cor}[thm]{Corollary}
\newtheorem{lem}[thm]{Lemma}
\newtheorem{prop}[thm]{Proposition}
\theoremstyle{definition}
\newtheorem{defn}[thm]{Definition}
\theoremstyle{remark}
\newtheorem{rem}[thm]{Remark}
\theoremstyle{remark}

\numberwithin{equation}{section}

\begin{document}

\title[Path decompositions and random times]{Enlargements of filtrations and path decompositions at non stopping times}%
\author{Ashkan Nikeghbali}
\address{American Institute of Mathematics
360 Portage Ave Palo Alto, CA 94306-2244 \\
and University of Rochester} \email{ashkan@aimath.org}
\date{\today}
\begin{abstract}
Az\'{e}ma associated with an honest time $L$ the supermartingale
$Z_{t}^{L}=\mathbb{P}\left[L>t|\mathcal{F}_{t}\right]$ and
established some of its important properties. This supermartingale
plays a central role in the general theory of stochastic processes
and in particular in the theory of progressive enlargements of
filtrations. In this paper, we shall give an additive
characterization for these supermartingales, which in turn will
naturally provide  many examples of enlargements of filtrations. We
combine this characterization with some arguments from both initial
and progressive enlargements of filtrations  to establish some path
decomposition results, closely related to or reminiscent of
Williams' path decomposition results. In particular, some of the
fragments of the paths in our decompositions end or start with a new
family of random times which are not stopping times, nor honest
times.
\end{abstract}
\maketitle

\section{Introduction}
Other than stopping times, the most commonly studied random times
occur as the ends of optional sets, or honest times. We shall denote
the class of such times by $L$. A very powerful technique for
studying such random times is that of progressive enlargement of
filtrations. The theory of progressive enlargements of filtrations
was introduced independently by Barlow (\cite{barlow}) and Yor
(\cite{gros}), and further developed by Jeulin and Yor
(\cite{jeulinyor,yorjeulin,jeulin,zurich}). It has many applications
in various parts of Probability Theory: path decompositions for some
diffusion processes (\cite{jeulin}, \cite{ashyordoob}), mathematical
models of default times and insider trading in mathematical finance
(\cite{elliotjeanbyor}), probabilistic inequalities
(\cite{jeulinyor},\cite{Ashkanbdg}), or new proofs of well known
results, such as Pitman's theorem (see \cite{zurich}, chapter XII).

Let
$\left(\Omega,\mathcal{F},\left(\mathcal{F}_{t}\right),\mathbb{P}\right)$
be a filtered probability space, satisfying the usual assumptions,
and $L$ the end of an
$\left( \mathcal{F}_{t}\right) $ optional set $%
\Gamma $, i.e:
\begin{equation*}
L=\sup \left\{ t:\left( t,\omega \right) \in \Gamma \right\}.
\end{equation*} The main idea is to consider the larger filtration
$$\mathcal{F}_{t}^{L}=\mathcal{G}_{t+},\;\;\mathrm{with}\;\;\quad\mathcal{G}_{t}\equiv\mathcal{F}_{t}\vee
\sigma\left\{L\wedge t,\right\},$$ which is the smallest right
continuous filtration which contains $\left(\mathcal{F}_{t}\right)$
and which makes $L$ a stopping time, and then to see how martingales
of the smaller filtration are changed when considered as stochastic
processes of the larger one. One process plays an essential role in
this theory, namely the supermartingale:
\begin{equation}
    Z_{t}^{L}=\mathbb{P}\left(L>t|\mathcal{F}_{t}\right),
\end{equation}
associated with $L$ by Az\'ema in \cite{azema}, and chosen to be
c\`{a}dl\`{a}g. An $\left(\mathcal{F}_{t}\right)$ local martingale
$\left(M_{t}\right)$, is a semimartingale in the larger filtration
$\left(\mathcal{F}_{t}^{L}\right)$ and decomposes as:
\begin{equation}\label{formuledegrossissiment}
M_{t}=\widetilde{M}_{t}+\int_{0}^{t\wedge L}\frac{d\langle M,Z^{L}\rangle_{s}}{Z_{s_{-}}^{L}}%
-\int_{L}^{t}\frac{d\langle M,Z^{L}\rangle_{s}}{1-Z_{s_{-}}^{L}%
},
\end{equation}
where $\left( \widetilde{M}_{t}\right) _{t\geq 0}$ denotes a $\left(
\left( \mathcal{F}_{t}^{L}\right) ,\mathbb{P}\right) $ local
martingale. One limitation of this formula is that it may not be
easy to compute the  supermartingale $Z^{L}$ for a given $L$; in
fact only a few examples are known (see \cite{jeulin, zurich}).
Hence, it would be useful to obtain a characterization result which
would help produce examples of honest times and their associated
supermartingales.

For simplicity, we make the following assumptions throughout this
paper, which we call the $\mathbf{(CA)}$ conditions:
\begin{enumerate}
\item all $\left( \mathcal{F}_{t}\right) $-martingales are \textbf{\underline{c}}ontinuous (e.g:
 the Brownian filtration).

\item the random time $L$ \textbf{\underline{a}}voids every $\left( \mathcal{F}_{t}\right) $%
-stopping time $T$, i.e. $\mathbb{P}\left[ L =T\right] =0$.\bigskip
\end{enumerate}
\begin{rem}
Under the conditions $\mathbf{(CA)}$, the optional and the
predictable sigma fields (with respect to
$\left(\mathcal{F}_{t}\right)$) are equal and the supermartingale
$(Z_{t}^{L})$ is continuous.
\end{rem}
One of the aims of this paper is to  characterize  the
supermartingales $(Z_{t}^{L})$. In \cite{ashyordoob}, we gave a
multiplicative characterization for the supermartingales
$(Z_{t}^{L})$, while here we shall adopt an additive approach
(Doob-Meyer decomposition). The paper is organized as follows:

In Section 2, we prove the characterization result for Az\'{e}ma's
supermartingales. To state it, we need to define a special class of
submartingales, whose definition goes back to Yor \cite{yorinegal},
and which was also studied in  \cite{Ashkanssmartrem} (under more
general conditions):
\begin{defn}\label{martreflechies}
Let $\left(X_{t}\right)$ be a positive local submartingale, which
decomposes as:
$$X_{t}=N_{t}+A_{t}.$$
We say that $\left(X_{t}\right)$ is of class $(\Sigma)$ if:
\begin{enumerate}
\item $\left(N_{t}\right)$ is a continuous local martingale, with $N_{0}=0$;
\item $\left(A_{t}\right)$ is a continuous increasing process, with $A_{0}=0$;
\item the measure $\left(dA_{t}\right)$ is carried by the set
$\left\{t:\;X_{t}=0\right\}$.
\end{enumerate}If additionally, $\left(X_{t}\right)$
is of class $(D)$, we shall say that $\left(X_{t}\right)$ is of
class $(\Sigma D)$.
\end{defn}Now,
consider the Doob-Meyer decomposition of $Z_{t}^{L}$:
$$Z_{t}^{L}\equiv1+\mu_{t}^{L}-A_{t}^{L}.$$We  prove that if
 $\left(Z_{t}\right)$ is a nonnegative supermartingale, with
$Z_{\infty}=0$, then, $Z$ may be represented as
$\mathbb{P}\left(L>t|\mathcal{F}_{t}\right)$, for some honest time
$L$ which avoids stopping times, if and only if
$\left(X_{t}\equiv1-Z_{t}\right)$ is a submartingale of the class
$(\Sigma )$, with the limit condition:
$$\lim_{t\rightarrow\infty}X_{t}=1.$$

\noindent  Section 3 contains our main results: we apply the results
of Section 2 to obtain decompositions analogous to  Williams' path
decomposition result for the supermartingale
$\left(Z_{t}^{L}\right)$. We also establish some path decomposition
results for certain classes of diffusion processes which play an
important role in applications. In particular, we shall see that the
pseudo-stopping times, introduced in \cite{AshkanYor}, play an
important role in path decompositions, exhibiting thus a new family
of random times enjoying nice properties with respect to path
decomposition.

\section{A characterization of Az\'{e}ma's supermartingale and applications}
\subsection{The characterization of Az\'{e}ma's supermartingale for honest times}
Az\'ema has studied in depth the  supermartingale
$Z_{t}^{L}=\mathbb{P}\left(L>t|\mathcal{F}_{t}\right)$ associated
with an honest time $L$ and has proved many interesting properties.
A classical example of such a random time, which has received much
attention in the literature (\cite{jeulin,zurich}, see
\cite{Ashkanbessel} for a one parameter extension),
is:$$L=\sup\left\{t\leq1:\; B_{t}=0\right\},$$ where as usual
$\left(B_{t}\right)$ denotes the standard Brownian Motion.

Let us briefly recall some results of Az\'{e}ma. We assume that the
conditions  $\mathbf{(CA)}$ hold. We  consider the Doob-Meyer
decomposition of $Z^{L}$:
\begin{equation}\label{doobmeyer}
Z_{t}^{L}=1+\mu _{t}^{L}-A_{t}^{L}
\end{equation}
The process $\left( A_{t}^{L}\right)$, which we shall sometimes
denote $\left( A_{t}\right)$ in the sequel, is the dual predictable
projection of the increasing process $\mathbf{1}_{\left\{ L\leq
t\right\} }$, and
\begin{equation*}
\mu_{t}^{L}=\mathbb{E}\left( A_{\infty }^{L}\mid
\mathcal{F}_{t}\right) -1
\end{equation*}
\begin{prop}[Az\'{e}ma \cite{azema}]\label{lemazema}
Let $L$ be the end of an optional set, or  an honest time (as was
discovered by M.T. Barlow \cite{barlow}, every honest time is the
end of some optional set); then
$$L=\sup\left\{t:\;Z_{t}=1\right\},$$and the measure $dA_{t}$ is
carried by the set $\left\{t:\;Z_{t}=1\right\}$. In particular, $A$
does not increase after $L$, i.e. $A_{L}=A_{\infty}$.
\end{prop}
\bigskip

To prove our main theorem, we shall need the following very useful
lemma, which  appears in the papers of Az\'ema, Meyer and Yor
\cite{azemameyeryor} and Az\'ema and Yor \cite{azemayorzero}.
\begin{lem}\label{martingalerefelchies}
Let $\left(X_{t}\right)$ be a submartingale of the class $(\Sigma
D)$ and let
$$L=\sup\left\{t:\;X_{t}=0\right\}.$$ Assume further that:
$$\mathbb{P}\left(X_{\infty}=0\right)=0.$$ Then:
\begin{equation}\label{martrel}
    X_{t}=\mathbb{E}\left(X_{\infty}\mathbf{1}_{\left\{L\leq
    t\right\}}|\mathcal{F}_{t}\right).
\end{equation}
\end{lem}
Now, we can state our characterization theorem.
\begin{thm} \label{charthm}
Let $\left(X_{t}\right)$ be a submartingale of the class $(\Sigma
D)$ satisfying: $\displaystyle\lim_{t\rightarrow\infty}X_{t}=1$. Let
$$L=\sup\left\{t:\;X_{t}=0\right\}.$$ Then $\left(X_{t}\right)$ is related to  Az\'{e}ma's supermartingale associated with $L$ in the
following way:
$$X_{t}=1-Z_{t}^{L}=\mathbb{P}\left(L\leq
t|\mathcal{F}_{t}\right).$$Consequently, if
 $\left(Z_{t}\right)$ is a nonnegative supermartingale, with
$Z_{0}=1$, then, $Z$ may be represented as
$\mathbf{P}\left(L>t|\mathcal{F}_{t}\right)$, for some honest time
$L$ which avoids stopping times, if and only if
$\left(X_{t}\equiv1-Z_{t}\right)$ is a submartingale of the class
$(\Sigma )$, with the limit condition:
$$\lim_{t\rightarrow\infty}X_{t}=1.$$
\end{thm}
\begin{proof}
This is an immediate application of Lemma \ref{martrel}, with
$X_{\infty}=1$ and Proposition \ref{lemazema}.
\end{proof}
\subsection{Some fundamental examples}
In the sequel, we give some explicit (yet generic) computations of
Az\'{e}ma's supermartingales for some honest times associated with
some very well known stochastic processes. These computations are
the first steps towards the path decompositions proved in the next
section.
\subsubsection{A Brownian example}
First, consider $\left(B_{t}\right)$, the standard Brownian Motion,
and let $T_{1}=\inf\left\{t\geq 0:\;B_{t}=1\right).$ Let
$\sigma=\sup\left\{t<T_{1}:\;B_{t}=0\right\}$. Then $B_{t\wedge
T_{1}}^{+}$ satisfies the conditions of Theorem \ref{charthm}, and
hence: $$\mathbb{P}\left(\sigma\leq
t|\mathcal{F}_{t}\right)=B_{t\wedge T_{1}}^{+}=\int_{0}^{t\wedge
T_{1}}\mathbf{1}_{B_{u}>0}\;dB_{u}+\dfrac{1}{2}\ell_{t\wedge
T_{1}},$$where $\left(\ell_{t}\right)$ is the local time of $B$ at
$0$. This example plays an important role in Williams' celebrated
path decomposition for the standard Brownian Motion on $[0,T_{1}]$.
This result is usually obtained by exploiting the strong Markov
property of the Brownian Motion. Our method allows us to get rid of
the Markov property, and to get similar formulae in the more general
context of continuous local martingales, as is shown in the next
paragraph.

One can also consider $T_{\pm1}=\inf\left\{t\geq
0:\;|B_{t}|=1\right\}$ and
$\tau=\sup\left\{t<T_{\pm1}:\;|B_{t}|=0\right\}$. $|B_{t\wedge
T_{\pm1}}|$ satisfies the conditions of Theorem \ref{charthm}, and
hence: $$\mathbb{P}\left(\tau\leq
t|\mathcal{F}_{t}\right)=|B_{t\wedge T_{\pm1}}|=\int_{0}^{t\wedge
T_{\pm1}}sgn\left(B_{u}\right)dB_{u}+\ell_{t\wedge T_{\pm1}}.$$
\subsubsection{Generalization to continuous local martingales}
More generally, consider a continuous local martingale
$\left(M_{t}\right)$ such that $M_{0}=0$ and
$<M>_{\infty}=\infty,\;\mathrm{a.s.}$; let $T_{1}=\inf\left\{t\geq
0:\;M_{t}=1\right\}$ and
$\sigma=\sup\left\{t<T_{1}:\;M_{t}=0\right\}$. Then, again, an
application of Theorem \ref{charthm} gives:
$$\mathbb{P}\left(\sigma\leq t|\mathcal{F}_{t}\right)=M_{t\wedge
T_{1}}^{+}=\int_{0}^{t\wedge
T_{1}}\mathbf{1}_{M_{u}>0}\;dM_{u}+\dfrac{1}{2}L_{t\wedge
T_{1}},$$where $\left(L_{t}\right)$ is the local time of $M$ at $0$.
\subsubsection{Recurrent diffusions}\label{paragraphrecdiff}
Let $\left(Y_{t}\right)$ be a real continuous recurrent diffusion
process, with $Y_{0}=0$. Then from the general theory of diffusion
processes, there exists a unique continuous and strictly increasing
function $s$, with $s\left(0\right)=0$, $\lim_{x\rightarrow
+\infty}s\left(x\right)=+\infty$, $\lim_{x\rightarrow
-\infty}s\left(x\right)=-\infty$, such that $s\left(Y_{t}\right)$ is
a continuous local martingale. Our aim is to establish some results
analogous to those established for the Brownian Motion and recurrent
continuous local martingales. Let
$$T_{1}\equiv\inf\left\{t\geq 0:\;Y_{t}=1\right).$$ Now, if we
define
$$X_{t}\equiv\dfrac{s\left(Y_{t\wedge
T_{1}}\right)^{+}}{s\left(1\right)},$$we easily note that $X$ is a
local submartingale of the class $(\Sigma)$ which satisfies the
hypotheses of Theorem \ref{charthm}. Consequently, if we note
$$\sigma=\sup\left\{t<T_{1}:\;Y_{t}=0\right\},$$we have:
$$\mathbb{P}\left(\sigma\leq
t|\mathcal{F}_{t}\right)=\dfrac{s\left(Y_{t\wedge
T_{1}}\right)^{+}}{s\left(1\right)}.$$
\subsubsection{Nonnegative continuous martingales which vanish at infinity}
Now let $\left(M_{t}\right)$ be a positive local martingale, such
that: $M_{0}=x,\; x>0$ and $\lim_{t\rightarrow\infty}M_{t}=0$. Then,
Tanaka's formula shows  that $\left(1-\dfrac{M_{t}}{y}\wedge
1\right)$, for $0\leq y \leq x$, is a local submartingale of the
class $(\Sigma)$ satisfying the assumptions of Theorem
\ref{charthm}, and hence with
$$g=\sup\left\{t:\; M_{t}=y\right\},$$we have:
$$\mathbb{P}\left(g> t|\mathcal{F}_{t}\right)=\dfrac{M_{t}}{y}\wedge
1=1+\dfrac{1}{y}\int_{0}^{t}\mathbf{1}_{\left(M_{u}<y\right)}\;dM_{u}-\dfrac{1}{2y}L_{t}^{y},$$where
$\left(L_{t}^{y}\right)$ is the local time of $M$ at $y$.
\subsubsection{Transient diffusions}\label{paragraphediffustrans}
As an illustration of the previous example, consider
$\left(R_{t}\right)$, a transient diffusion with values in
$\left[0,\infty\right)$, which has $\left\{0\right\}$ as entrance
boundary. Let $s$ be a scale function for $R$, which we can choose
such that: $$s\left(0\right)=-\infty, \text{ and }
s\left(\infty\right)=0.$$ Then, under the law $\mathbb{P}_{x}$, for
any $x>0$, the local martingale
$\left(M_{t}=-s\left(R_{t}\right)\right)$ satisfies the conditions
of the previous example and for $0\leq x \leq y$, we
have:
\begin{equation}\label{pitmanyorbessel}
    \mathbb{P}_{x}\left(g_{y}>
t|\mathcal{F}_{t}\right)=\dfrac{s\left(R_{t}\right)}{s\left(y\right)}\wedge
1=1+\dfrac{1}{s\left(y\right)}\int_{0}^{t}\mathbf{1}_{\left(R_{u}>y\right)}\;d\left(s\left(R_{u}\right)\right)+\dfrac{1}{2s\left(y\right)}L_{t}^{s\left(y\right)},
\end{equation}
where $\left(L_{t}^{s\left(y\right)}\right)$ is the local time of
$s\left(R\right)$ at $s\left(y\right)$, and where
$$g_{y}=\sup\left\{t:\; R_{t}=y\right\}.$$
Formula (\ref{pitmanyorbessel}) was the key point in the derivation
of the distribution of $g_{y}$ in \cite{pitmanyor}, Theorem 6.1,
p.326.
\section{Path decompositions}
In this section, inspired by Williams' path decompositions for the
standard Brownian Motion and for transient diffusions given their
minima, we establish  path decompositions for Az\'{e}ma's
supermartingales and some families of recurrent and transient
diffusions. What follows is similar in spirit to what we have done
in \cite{ashyordoob}, but in an additive setting, and of course the
results are different. It is also an opportunity to show that the
techniques of progressive and initial enlargements of filtrations
can be combined to prove, quite shortly, some non trivial path
decomposition results.

Let us recall briefly the random times introduced by D. Williams to
study the paths of a standard Brownian Motion $B$:
\begin{equation*}
T_{1}=\inf \left\{ t:\text{ }B_{t}=1\right\} ,\;\sigma =\sup \left\{
t<T_{1}:\text{ }B_{t}=0\right\} ;
\end{equation*}%
and
\begin{equation*}
\rho =\sup \left\{ u<\sigma :\text{ }B_{u}=S_{u}\right\} ,\text{ where }%
S_{u}=\sup_{s\leq u}B_{s}.
\end{equation*}D. Williams (\cite{williams}) discovered the remarkable fact that although
$\rho$ is not a stopping time, it  nevertheless satisfies the
optional stopping theorem, i.e. for every bounded martingale
$\left(M_{t}\right)$ of the filtration
$\left(\mathcal{F}_{t}\right)$, we have:
$$\mathbb{E}M_{\rho}=\mathbb{E}M_{\infty}.$$In \cite{AshkanYor}, we
have called such random times pseudo-stopping times and we have
characterized them.  Before stating and proving our results, we
shall first recall in the next subsection some standard facts about
pseudo-stopping times and multiple enlargements of filtrations.
\subsection{Basic facts about pseudo-stopping times and double enlargements of filtrations}
\subsubsection{Pseudo-stopping times}\label{secpta} In
\cite{AshkanYor}, following D. Williams, we have proposed the
following generalization of stopping times:
\begin{defn}[\cite{AshkanYor}]
Let $\rho:\;(\Omega,\mathcal{F})\rightarrow\mathbf{R}_{+}$ be a
random time; $\rho$ is called a pseudo-stopping time if for every
bounded $\left(\mathcal{F}_{t}\right)$ martingale
$\left(M_{t}\right)$ we have:
$$\mathbf{E}\left(M_{\rho}\right)=\mathbf{E}\left(M_{0}\right).$$
\end{defn}David Williams (\cite{williams}) gave the first example of such a random
time and the following systematic construction is established in
\cite{AshkanYor}:
\begin{prop}[\cite{AshkanYor}]\label{ptaconstruction}
Let $L$ be an honest time. Then, under the conditions \textbf{(CA)},
$$\rho\equiv \sup\left\{t<L:\;Z_{t}^{L}=\inf_{u\leq
L}Z_{u}^{L}\right\},$$is a pseudo-stopping time, with
$$Z_{t}^{\rho}\equiv
\mathbf{P}\left(\rho>t\mid\mathcal{F}_{t}\right)=\inf_{u\leq
t}Z_{u}^{L},$$and $Z_{\rho}^{\rho}$ follows the uniform distribution
on $(0,1)$.
\end{prop}The following property, also proved in \cite{AshkanYor}, is
essential in studying path decompositions:
\begin{prop}[\cite{AshkanYor}]\label{regenrative}
Let $\rho$ be a pseudo-stopping time and let $M_{t}$ be an
$\left(\mathcal{F}_{t}\right)$ local martingale. Then
$\left(M_{t\wedge \rho}\right)$ is an
$\left(\mathcal{F}_{t}^{\rho}\right)$ local martingale.
\end{prop}
To conclude, let us illustrate Proposition \ref{ptaconstruction}
with an example. Let $Y$ be a recurrent diffusion; with the
notations and assumptions of paragraph \ref{paragraphrecdiff},
$$\rho\equiv\sup\left\{t<\sigma:\;Y_{t}=\max_{u\leq\sigma}Y_{u}\right\},$$is
a pseudo-stopping time.
\subsubsection{Double enlargements of filtrations}
We recall some lesser known results of Jeulin (\cite{jeulin}) about
successive progressive enlargements of filtrations. The reader can
also refer to \cite{delmaismey} for a more recent exposition
(summary) of these facts.
\begin{prop}[Jeulin \cite{jeulin}]\label{resultdoublegrosjeulin}
Let $L$ be an honest time for the filtration
$\left(\mathcal{F}_{t}\right)$ and let $\rho$ be an honest time for
$\left(\mathcal{F}_{t}^{L}\right)$, and define
$\left(\mathcal{F}_{t}^{L,\rho}\right)$ the filtration obtained by
enlarging progressively $\left(\mathcal{F}_{t}^{L}\right)$ with
$\rho$. Then, any $\left(\mathcal{F}_{t}\right)$ local martingale
$\left(M_{t}\right)$ is a semimartingale in
$\left(\mathcal{F}_{t}^{L,\rho}\right)$ and decomposes as:
\begin{equation}\label{doublegross}
    M_{t}=\widetilde{M}_{t}+\int_{0}^{t\wedge\rho}\dfrac{d\langle
    M,Z^{\rho}\rangle_{u}}{Z_{u-}^{\rho}}+\int_{\rho}^{t\wedge L}\dfrac{d\langle
    M,Z^{L}-Z^{\rho}\rangle_{u}}{Z_{u-}^{L}-Z_{u-}^{\rho}}-\int_{L}^{t}\dfrac{d\langle
    M,Z^{L}\rangle_{u}}{1-Z_{u-}^{L}},
\end{equation}where $\left(\widetilde{M}_{t}\right)$ is an
$\left(\mathcal{F}_{t}^{L,\rho}\right)$ local martingale.
\end{prop}
Honest times enjoy the remarkable property that every
$\left(\mathcal{F}_{t}\right)$ semimartingale is an
$\left(\mathcal{F}_{t}^{L}\right)$ semimartingale, or in the jargon
of the theory of enlargements of filtrations, the pair of
filtrations
$\left(\left(\mathcal{F}_{t}\right),\left(\mathcal{F}_{t}^{L}\right)\right)$
satisfy the $\left(\mathcal{H}'\right)$ hypothesis. The previous
proposition shows that there might be non-honest times which enjoy
this property; indeed, the pseudo-stopping times defined in
Proposition \ref{ptaconstruction} have this property:
\begin{cor}
Let us consider the pseudo-stopping time defined in Proposition
\ref{ptaconstruction}. Then, every $\left(\mathcal{F}_{t}\right)$
semimartingale is an $\left(\mathcal{F}_{t}^{\rho}\right)$
semimartingale, or in other words, the pair of filtrations
$\left(\left(\mathcal{F}_{t}\right),\left(\mathcal{F}_{t}^{\rho}\right)\right)$
satisfy the $\left(\mathcal{H}'\right)$ hypothesis.
\end{cor}
\begin{proof}
It suffices to prove that every
$\left(\left(\mathcal{F}_{t}\right)\right)$ local martingale
$\left(M_{t}\right)$ is an $\left(\mathcal{F}_{t}^{\rho}\right)$
semimartingale. From Proposition \ref{resultdoublegrosjeulin}, every
$\left(\left(\mathcal{F}_{t}\right)\right)$ local martingale is an
$\left(\mathcal{F}_{t}^{L,\rho}\right)$ semimartingale, and since
$\mathcal{F}_{t}^{\rho}\subset\mathcal{F}_{t}^{L,\rho}$ and
$\left(M_{t}\right)$ is $\left(\mathcal{F}_{t}^{\rho}\right)$
adapted, it follows from a well known result of Stricker (see
\cite{protter}) that $\left(M_{t}\right)$ is also an
$\left(\mathcal{F}_{t}^{\rho}\right)$ semimartingale.
\end{proof}
We shall also need another result of Jeulin which certainly deserves
to be better known: the problem of initial enlargement with
$A_{\infty}^{L}$.
\begin{prop}\label{grosinitialavecainfalajeulin}
Let $T$ be a totally inaccessible stopping time, such that
$\mathbb{P}\left(T>0\right)=1$ and let $\left(A_{t}\right)$ be the
$\left(\mathcal{F}_{t}\right)$ dual predictable projection of
$\left(\mathbf{1}_{T\leq t}\right)$. Then the following hold:
\begin{enumerate}
\item $A$ is continuous and
$T=\inf\left\{t:\;A_{t}=A_{T}\right\}=\sup\left\{t:\;A_{t}=A_{T}\right\}$
(Az\'{e}ma \cite{azema});
\item define
$\mathcal{G}_{t}\equiv\bigcap_{\varepsilon>0}\left(\mathcal{F}_{t+\varepsilon}\vee\sigma\left(A_{T}\right)\right)$;
then every continuous $\left(\mathcal{F}_{t}\right)$ martingale is a
$\left(\mathcal{G}_{t}\right)$
 martingale (Jeulin \cite{jeulin}).
\end{enumerate}
\end{prop}
\begin{rem}
In fact, every $\left(\mathcal{F}_{t}\right)$ martingale, which does
not jump at $T$, is a $\left(\mathcal{G}_{t}\right)$
 martingale.
\end{rem}
\begin{cor}\label{remarquesurletrhocommetpsaret}
\begin{enumerate}
\item If $L$ is an honest time which avoids stopping times, then every $\left(\mathcal{F}_{t}\right)$ martingale is an
$\left(\mathcal{F}_{t}^{L,\sigma\left(A_{\infty}^{L}\right)}\right)$
semimartingale and the decomposition formula is the same as the
decomposition formula (\ref{formuledegrossissiment}):
\begin{equation*}M_{t}=\widetilde{M}_{t}+\int_{0}^{t\wedge L}\frac{d\langle M,Z^{L}\rangle_{s}}{Z_{s_{-}}^{L}}%
+\int_{L}^{t}\frac{d\langle M,Z^{L}\rangle_{s}}{1-Z_{s_{-}}^{L}%
}.
\end{equation*}
\item Similarly, under the assumptions $(\mathbf{CA})$, the
pseudo-stopping times of Proposition \ref{ptaconstruction} are
inaccessible stopping times for the filtration
$\left(\mathcal{F}_{t}^{L,\rho}\right)$ and
$\left(\log\left(\frac{1}{Z_{t\wedge\rho}^{\rho}}\right)\right)$ is
the $\left(\mathcal{F}_{t}^{L,\rho}\right)$ dual predictable
projection of $\left(\mathbf{1}_{\rho\leq t}\right)$. Here again,
every $\left(\mathcal{F}_{t}\right)$ martingale is an
$\left(\mathcal{F}_{t}^{L,\sigma\left(Z_{\rho}^{\rho}\right)}\right)$
semimartingale whose decomposition is given by (\ref{doublegross}):
\begin{equation}\label{ptadecomp}
M_{t}=\widetilde{M}_{t}+\int_{\rho}^{t\wedge L}\dfrac{d\langle
    M,Z^{L}\rangle_{u}}{Z_{\rho}^{L}-Z_{u}^{\rho}}-\int_{L}^{t}\dfrac{d\langle
    M,Z^{L}\rangle_{u}}{1-Z_{u}^{L}},
\end{equation}where $\left(\widetilde{M}_{t}\right)$ is an
$\left(\mathcal{F}_{t}^{L,\rho}\right)$ and
$\left(\mathcal{F}_{t}^{L,\sigma\left(Z_{\rho}^{\rho}\right)}\right)$
local martingale, every continuous
$\left(\mathcal{F}_{t}^{L,\rho}\right)$ martingale being an
$\left(\mathcal{F}_{t}^{L,\sigma\left(Z_{\rho}^{\rho}\right)}\right)$
martingale.
\end{enumerate}
\end{cor}
\begin{proof}
$(1).$ This first point follows from the fact that $L$ is a totally
inaccessible stopping time for $\left(\mathcal{F}_{t}^{L}\right)$
(see \cite{yorjeulin}), and  Proposition
\ref{grosinitialavecainfalajeulin} can be applied with $A_{T}\equiv
A_{\infty}^{L}$.

$(2).$ First, we note from Proposition \ref{ptaconstruction} that
$\left(Z_{t}^{\rho}\right)$ is a continuous and decreasing process
($Z_{t}^{\rho}=1-A_{t}^{\rho}$). Moreover, from a result of Jeulin
and Yor (\cite{yorjeulin}), $\mathbf{1}_{\left(\rho\leq
t\right)}-\int_{0}^{t\wedge\rho}\frac{A_{u}^{\rho}}{Z_{u}^{\rho}}=\mathbf{1}_{\left(\rho\leq
t\right)}-\log\left(\frac{1}{Z_{t\wedge\rho}^{\rho}}\right)$ is an
$\left(\mathcal{F}_{t}^{\rho}\right)$ martingale. It remains a
martingale in $\left(\mathcal{F}_{t}^{L,\rho}\right)$, since it is
of finite variation and $\rho<L$. Consequently,
$\left(\log\left(\frac{1}{Z_{t\wedge\rho}^{\rho}}\right)\right)$ is
also the $\left(\mathcal{F}_{t}^{L,\rho}\right)$ dual predictable
projection of $\left(\mathbf{1}_{\rho\leq t}\right)$ and the
announced results now easily follow from Propositions
\ref{resultdoublegrosjeulin} and \ref{grosinitialavecainfalajeulin}.
\end{proof}

\subsection{An analogue of Williams' path decomposition theorem for $\mathbf{\mathbb{P}\left(L\leq
t|\mathcal{F}_{t}\right)}$ } We are going to use techniques from
both stochastic calculus and the general theory of stochastic
processes (the Dubins-Schwarz theorem and the decomposition formula
in the larger filtration) to generalize some fragments of Williams'
path decomposition for the standard Brownian to more general
processes, namely the submartingale $\left(\mathbb{P}\left(L\leq
t|\mathcal{F}_{t}\right)\right)$, associated with an honest time
$L$, under the conditions $\mathbf{(CA)}$.

\begin{thm}\label{decogenraalawill}
Let $$X_{t}=N_{t}+A_{t},$$be a submartingale of the class $(\Sigma)$
satisfying \begin{equation} \label{condi1}
\lim_{t\rightarrow\infty}X_{t}=1. \end{equation} and let
$L=\sup\left\{t:\;X_{t}=0\right\}.$ Recall (Theorem \ref{charthm})
that
$$X_{t}=1-Z_{t}^{L}=\mathbf{P}\left(L\leq
t|\mathcal{F}_{t}\right).$$  Let us also define the random
time:$$\rho=\sup \left\{ t<L:\text{ }X_{t}=S_{L}\right\},$$where
$$S_{t}=\sup_{u\leq
t}X_{u}.$$Then:
\begin{enumerate}
\item the process
$\left(X_{t}\right)$, stopped at $\rho$ is, up to the time change
$\left(\langle N \rangle_{t}\right)$, a reflected Brownian Motion
started from $0$, stopped when it first hits $1$.
\item The random time $\rho$ is a pseudo-stopping time and $$\mathbb{P}\left(\rho>t|\mathcal{F}_{t}\right)=1-S_{t}.$$ Moreover, $X_{\rho}$ is uniformly distributed on $(0,1)$, and
conditionally on $X_{\rho}=m$, $\left(X_{t}\right)$ is up to the
time change $\left(\langle N \rangle_{t}\right)$, a reflected
Brownian Motion started from $0$ and stopped when it first hits $m$.
\item Conditionally on $\mathcal{F}_{L}$, the process
$\left(X_{L+t}\right)$ is, up to the time change \linebreak
$\left(\langle N \rangle_{L+t}-\langle N \rangle_{L}\right)$, a
Bessel process of dimension $3$, started from $0$, and stopped when
it first hits $1$.
\end{enumerate}
\end{thm}
\begin{proof}
$\left(1\right)$. First, from  Skorokhod's reflection lemma (see
\cite{revuzyor} or \cite{Ashkanssmartrem}), we have:
$$A_{t}=\sup_{u\leq t}\left(-N_{u}\right).$$ Moreover, there exists
a Brownian Motion $\left(\beta_{u}\right)$ such that:
$$N_{t}=\beta_{\langle N
\rangle_{t}}.$$ Hence, up to the time change $\left(\langle N
\rangle_{t}\right)$, $\left(X_{t}\right)$ has the same decomposition
as the absolute value of a Brownian Motion (this is immediate from
Tanaka's formula). Thus it is a time changed reflected Brownian
Motion.

$\left(2\right)$. The first point follows immediately from
Proposition \ref{ptaconstruction}: indeed, $\rho$ is a
pseudo-stopping time and $X_{\rho}$ is equal to $Z_{\rho}^{\rho}$,
which is uniformly distributed. The second point follows from a
combination of  Proposition \ref{resultdoublegrosjeulin}
 and Corollary
\ref{remarquesurletrhocommetpsaret}. Indeed, in the filtration
$\left(\mathcal{F}_{t}^{L,\sigma\left(X_{\rho}\right)}\right)$,
obtained by initially enlarging the filtration
$\left(\mathcal{F}_{t}^{L}\right)$ with
$\log\left(\frac{1}{Z_{t\wedge\rho}^{\rho}}\right)=\log\left(\frac{1}{1-X_{\rho}}\right)$,
we have:
$$X_{t\wedge\rho}=N_{t\wedge\rho}+A_{t\wedge\rho}.$$

$\left(3\right)$. We first note that, since $X_{L}=0$, we have
$N_{L}=-A_{L}$, and consequently, $$X_{L+t}=N_{L+t}-N_{L}.$$ Now,
using the fact that $X_{t}=1-Z_{t}^{L}=\mathbf{P}\left(L\leq
t|\mathcal{F}_{t}\right)$ the decomposition formula
(\ref{formuledegrossissiment}) yields:
$$X_{L+t}=N_{L+t}-N_{L}=\widetilde{N}_{t}+\int_{0}^{t}\frac{d\langle N
\rangle_{L+u}}{X_{L+u}},$$ where $\widetilde{N}$ is an
$\left(\mathcal{F}_{t}^{L}\right)$ martingale. Now, the result
follows from the fact that the Bessel process of dimension $3$
$\left(R_{t}\right)$ can be characterized as the unique solution to
the stochastic differential equation:
$$dR_{t}=dB_{t}+\dfrac{dt}{R_{t}},$$ where
$\left(B_{t}\right)$ is a one dimensional Brownian Motion.
\end{proof}
As an illustration of the above theorem, let us consider
$$X_{t}\equiv\alpha B_{t}^{+}+\beta B_{t}^{-},$$where $B$ is the
standard Brownian Motion and $\alpha>0,\beta>0$. Let
$T_{1}=\inf\left\{t:\;X_{t}=1\right\}$. Then, it is easy to check
that $\left(X_{t\wedge T_{1}}\right)$ satisfies the assumptions of
the Theorem \ref{decogenraalawill} with the time change $$\langle N
\rangle_{t}=\alpha^{2}\int_{0}^{t}\mathbf{1}_{\left(B_{u}>0\right)}\;du+\beta^{2}\int_{0}^{t}\mathbf{1}_{\left(B_{u}\leq0\right)}\;du.$$
\subsection{Path decompositions for some recurrent diffusions}
D. Williams's path decomposition also admits a  generalization to
the  wider class of recurrent diffusions $\left(Y_{t}\right)$,
satisfying the stochastic differential equation:
\begin{equation}\label{equationrecurrence}
Y_{t}=B_{t}+\int_{0}^{t}b\left(Y_{u}\right)du,
\end{equation}where $\left(B_{t}\right)$ is the standard Brownian
Motion, and $b$ is a Borel integrable function which allows
existence and uniqueness for equation (\ref{equationrecurrence}). We
note $\mathcal{L}$ the infinitesimal generator of this diffusion:
$$\mathcal{L}=\frac{1}{2}\dfrac{d^{2}}{dx^{2}}+b\left(x\right)\dfrac{d}{dx}.$$Let $T_{1}\equiv
\inf\left\{t:\;Y_{t}=1\right)$, and denote by $s$ the scale function
of $Y$, which is strictly increasing and which vanishes at zero,
i.e:
$$s\left(z\right)=\int_{0}^{z}\exp\left(-2\widehat{b}\left(y\right)\right)dy,$$where$$\widehat{b}\left(y\right)=\int_{0}^{y}b\left(u\right)du.$$ From the results
of paragraph (\ref{paragraphrecdiff}), if we define
$$\sigma=\sup\left\{t<T_{1}:\;Y_{t}=0\right\},$$we have, with  $\left(L_{t}\right)$
the local time at $0$ of the local martingale $s\left(Y\right)$:
\begin{eqnarray*}
  \mathbb{P}\left(\sigma\leq
t|\mathcal{F}_{t}\right) &=& \dfrac{s\left(Y_{t\wedge
T_{1}}\right)^{+}}{s\left(1\right)} \\
   &=& \dfrac{1}{s\left(1\right)}\int_{0}^{t\wedge
T_{1}}\mathbf{1}_{\left(Y_{u}>0\right)}s'\left(Y_{u}\right)dB_{u}+\dfrac{1}{2s\left(1\right)}L_{t\wedge
T_{1}},
\end{eqnarray*}where the last equality is obtained by an application
of Tanaka's formula. Moreover, from Proposition
(\ref{ptaconstruction}),
\begin{equation}\label{ptarecdiff}
    \rho\equiv\sup\left\{t<\sigma:\;Y_{t}=\max_{u\leq\sigma}Y_{u}\right\},
\end{equation}is a pseudo-stopping time.
\begin{prop}\label{decodifrec}
Let $\left(Y_{t}\right)$ be a diffusion process satisfying equation
(\ref{equationrecurrence}). Define:
$$\overline{Y}_{t}=\max_{u\leq t}Y_{u}.$$
Then:
\begin{itemize}
\item The process $\left(Y_{\sigma+t},\;t\leq T_{1}-\sigma\right)$
is an $\left(\mathcal{F}_{\sigma+t}\right)$ diffusion, starting from
$0$, considered up to the first time it hits $1$, and is independent
of $\mathcal{F}_{\sigma}$. Its infinitesimal generator is given by:
$$\mathcal{L}=\dfrac{1}{2}\dfrac{d^{2}}{dx^{2}}+\left(b\left(x\right)+\dfrac{s'\left(x\right)}{s\left(x\right)}\right)\dfrac{d}{dx}.$$
\item The random time $\rho$ is a pseudo-stopping time and satisfies:
$$\mathbb{P}\left(\rho>t|\mathcal{F}_{t}\right)=1-\dfrac{s\left(\overline{Y}_{t\wedge T_{1}}\right)^{+}}{s\left(1\right)}.$$
Moreover, $Y_{\rho}=\overline{Y}_{\sigma}$ follows the same law as
$s^{-1}\left(s\left(1\right)U\right)$, where $U$ is a random
variable following the uniform law on $(0,1)$, and is independent of
the whole process $\left(Y_{\sigma+t},\;t\leq T_{1}-\sigma\right)$.
\item \textbf{Conditionally on} $\mathbf{Y_{\rho}=m}$,
\begin{enumerate}
\item the process $\left(Y_{t};\;t\leq \rho\right)$ is a diffusion process, considered up to $T_{m}$, the first time when it
hits $m$, with the same infinitesimal generator as $Y$.
\item the process $\left(Y_{\rho+t};\;t\leq
\sigma-\rho\right)$ is a $\left(\mathcal{F}_{\rho+t}\right)$
diffusion process, started from $m$, considered up to $T_{0}$, the
first time when it hits $0$, and is independent of
$\left(Y_{t};\;t\leq \rho\right)$; its infinitesimal generator is
given by:
$$\frac{1}{2}\dfrac{d^{2}}{dx^{2}}+\left(b\left(x\right)+\mathbf{1}_{\left(x>0\right)}\dfrac{s'\left(x\right)}{s\left(x\right)-s\left(m\right)}\right)\dfrac{d}{dx}.$$
\end{enumerate}
\end{itemize}
\end{prop}
\begin{proof}
The proof is based on enlargements arguments. First, we note that
from Proposition \ref{ptaconstruction} $s\left(Y_{\rho}\right)$ is
distributed as $s\left(1\right)U$, where $U$ follows the uniform law
on $(0,1)$.

Now, let us study the path of $Y$ on $[\sigma,T_{1}]$. From formula
(\ref{formuledegrossissiment}), the Brownian Motion $B$ is a
semimartingale in the filtration
$\left(\mathcal{F}_{t}^{\sigma}\right)$ and decomposes as:
$$B_{t}=\widetilde{B}_{t}+\int_{0}^{t\wedge\sigma}\dfrac{d<B,Z^{\sigma}>_{u}}{Z_{u}^{\sigma}}%
+\int_{\sigma}^{t\wedge T_{1}}\frac{d<B,1-Z^{\sigma}>_{u}}{1-Z_{u}^{\sigma}%
},$$where $\widetilde{B}$ is a
$\left(\mathcal{F}_{t}^{\sigma}\right)$ Brownian Motion (indeed it
is a continuous local martingale with bracket $t$). Consequently,
the diffusion $Y$, which is an$\left(\mathcal{F}_{t}\right)$
semimartingale, remains a semimartingale in
$\left(\mathcal{F}_{t}^{\sigma}\right)$ and its decomposition is
given by:
$$Y_{t}=\widetilde{B}_{t}+\int_{0}^{t}b\left(Y_{u}\right)du-\int_{0}^{t\wedge\sigma}\mathbf{1}_{Y_{u}>0}\dfrac{s'\left(Y_{u}\right)}{s\left(1\right)-s\left(Y_{u}\right)}du+\int_{\sigma}^{t\wedge T_{1}}\dfrac{s'\left(Y_{u}\right)}{s\left(Y_{u}\right)}du.$$
Now, considering $Y_{\sigma+t}-Y_{\sigma}=Y_{\sigma+t}$, for $t\leq
T_{1}-\sigma$, we obtain:
$$Y_{\sigma+t}=\left(\widetilde{B}_{\sigma+t}-\widetilde{B}_{\sigma}\right)+\int_{\sigma}^{\sigma+t}b\left(Y_{u}\right)du+\int_{\sigma}^{(\sigma+t)\wedge T_{1}}\dfrac{s'\left(Y_{u}\right)}{s\left(Y_{u}\right)}du.$$
Now, using the fact that $\sigma$ is a stopping time for the
filtration $\left(\mathcal{F}_{t}^{\sigma}\right)$, we have that
$\left(\widetilde{B}_{\sigma+t}-\widetilde{B}_{\sigma}\right)$,
which we note $\left(\widehat{B}_{t}\right)$, is a Brownian Motion,
which is independent of
$\mathcal{F}_{\sigma}^{\sigma}\supset\mathcal{F}_{\sigma}$.
Consequently, for $t\leq T_{1}-\sigma$, we have:
$$Y_{\sigma+t}=\widehat{B}_{t}+\int_{0}^{t}b\left(Y_{\sigma+u}\right)du+\int_{0}^{t\wedge
\left(T_{1}-\sigma\right)}\dfrac{s'\left(Y_{\sigma+u}\right)}{s\left(Y_{\sigma+u}\right)}du,$$and
the result announced for the path on $[\sigma,T_{1}]$ follows now
easily.

Now, let us consider the path of $Y$ on $[0,\rho]$, and
$[\rho,\sigma]$. For this, we enlarge initially the filtration
$\left(\mathcal{F}_{t}^{\sigma}\right)$ with the variable
$Y_{\rho}$: according to Proposition \ref{resultdoublegrosjeulin}
and Corollary \ref{remarquesurletrhocommetpsaret}, for
$t\leq\sigma$, $B$ decomposes in
$\left(\mathcal{F}_{t}^{\sigma,\sigma\left(Y_{\rho}\right)}\right)$,
which we note $\left(\mathcal{F}_{t}^{\sigma,Y_{\rho}}\right)$ for
notational convenience, as:
$$B_{t}=\widetilde{B}_{t}-\int_{\rho}^{t\wedge\sigma}\mathbf{1}_{\left(Y_{u}>0\right)}\dfrac{s'\left(Y_{u}\right)}{s\left(Y_{\rho}\right)-s\left(Y_{u}\right)}du,$$
where $\widetilde{B}$ is an
$\left(\mathcal{F}_{t}^{\sigma,Y_{\rho}}\right)$ Brownian Motion
which is independent of $Y_{\rho}$. Consequently, for $t\leq\rho$,
$Y$ decomposes in $\left(\mathcal{F}_{t}^{\sigma,Y_{\rho}}\right)$
as:
\begin{equation}\label{eqa1}
    Y_{t}=\widetilde{B}_{t}+\int_{0}^{t}b\left(Y_{u}\right)du,\;\mathrm{for
    }\;t\leq\rho,
\end{equation}and for $t\leq(\sigma-\rho)$, we have:
\begin{equation}\label{eqa2}
 Y_{\rho+t}=Y_{\rho}+\left(\widetilde{B}_{\rho+t}-\widetilde{B}_{\rho}\right)+\int_{\rho}^{t\wedge(\sigma-\rho)}b\left(Y_{u}\right)du-\int_{\rho}^{t\wedge(\sigma-\rho)}\dfrac{s'\left(Y_{u}\right)}{s\left(Y_{\rho}\right)-s\left(Y_{u}\right)}du.
\end{equation}Now,
$$\widehat{B}_{t}\equiv\widetilde{B}_{\rho+t}-\widetilde{B}_{\rho}$$
is again a standard Brownian Motion, independent of $Y_{\rho}$, and
hence, conditionally on $Y_{\rho}=m$,  the process
$\left(Y_{\rho+t};\;t\leq \sigma-\rho\right)$ satisfies:
$$Y_{\rho+t}=m+\widehat{B}_{t}+\int_{0}^{t\wedge(\sigma-\rho)}b\left(Y_{\rho+u}\right)du+\int_{0}^{t\wedge(\sigma-\rho)}\dfrac{s'\left(Y_{\rho+u}\right)}{s\left(Y_{\rho+u}\right)-s\left(m\right)}du.$$
The statement of the Proposition now follows from the last equality
and equation (\ref{eqa1}).
\end{proof}
\begin{rem}
When $b\equiv0$, we have $s\left(x\right)=x$, and we recover D.
Williams'path decomposition for the standard Brownian Motion.
\end{rem}
\subsection{Path decompositions for some transient diffusions}
Now, we consider a special subfamily of the transient diffusions of
paragraph \ref{paragraphediffustrans} which play an important role
in the extension of Pitman's theorem (see \cite{zurich}, p.46). More
precisely, let $\left(R_{t}\right)$ be any transient diffusion which
takes its values in $\left(0,\infty\right)$, and satisfies:
\begin{equation}\label{edsdiftrans}
    R_{t}=x+B_{t}+\int_{0}^{t}c\left(R_{u}\right)du,\;,x>0\;t\geq0,
\end{equation}where $c:\mathbb{R}_{+}\rightarrow\mathbb{R}$ allows
uniqueness in law for this equation. Noting
$T_{0}=\inf\left\{t:\;R_{t}=0\right\}$, we assume that
$\mathbb{P}_{x}\left(T_{0}<\infty\right)=0$, so that a scale
function $s$ of $R$ may be chosen to satisfy:
$$s\left(0\right)=-\infty;s\left(\infty\right)=0;\;\dfrac{1}{2}s''+cs'=0.$$
We keep the notation of paragraph \ref{paragraphediffustrans}: for
$0\leq x \leq y$, and $$g_{y}=\sup\left\{t:\; R_{t}=y\right\},$$we
have:
\begin{eqnarray*}
  \mathbb{P}_{x}\left(g_{y}>
t|\mathcal{F}_{t}\right) &=&
\dfrac{s\left(R_{t}\right)}{s\left(y\right)}\wedge
1 \\
   &=& 1+\dfrac{1}{s\left(y\right)}\int_{0}^{t}\mathbf{1}_{R_{u}>y}\;s'\left(R_{u}\right)dB_{u}+\dfrac{1}{2s\left(y\right)}L_{t}^{s\left(y\right)},
\end{eqnarray*}where
$L_{t}^{s\left(y\right)}$ is the local time of $s\left(R\right)$ at
$s\left(y\right)$.

From Proposition \ref{ptaconstruction}, the random time:
$$\rho=\sup\left\{t<g_{y}:\;R_{t}=\sup_{u\leq g_{y}}R_{u}\right\},$$is a
pseudo-stopping time and
$\mathbb{P}\left(\rho>t|\mathcal{F}_{t}\right)=\dfrac{s\left(\sup_{u\leq
t}R_{u}\right)}{s\left(y\right)}\wedge 1$. Now, likewise Proposition
\ref{decodifrec}, the following path decomposition holds for the
diffusion $R$:
\begin{prop}
Let $\left(R_{t}\right)$ be a diffusion process satisfying equation
(\ref{edsdiftrans}). Then:
\begin{itemize}
\item The process $\left(R_{g_{y}+t},\;t\geq0\right)$
is an $\left(\mathcal{F}_{g_{y}+t}\right)$ diffusion, starting from
$y$,  and is independent of $\mathcal{F}_{g_{y}}$. Its infinitesimal
generator is given by:
$$\mathcal{L}=\dfrac{1}{2}\dfrac{d^{2}}{dx^{2}}+\left(c\left(x\right)+\dfrac{s'\left(x\right)}{s\left(x\right)-s\left(y\right)}\right)\dfrac{d}{dx}.$$
\item The random time $\rho$ is a pseudo-stopping time and  $R_{\rho}=\overline{R}_{g_{y}}$ follows the same law as
$s^{-1}\left(s\left(y\right)U\right)$, where $U$ follows the uniform
law on $(0,1)$, and is independent of the whole process
$\left(R_{g_{y}+t},\;t\geq0\right)$.
\item \textbf{Conditionally on} $\mathbf{R_{\rho}=m}$,
\begin{enumerate}
\item the process $\left(R_{t};\;t\leq \rho\right)$ is a diffusion process, considered up to $T_{m}$, the first time when it
hits $m$, with the same infinitesimal generator as $R$.
\item the process $\left(R_{\rho+t};\;t\leq
g_{y}-\rho\right)$ is a $\left(\mathcal{F}_{\rho+t}\right)$
diffusion process, started from $m$, considered up to $T_{y}$, the
first time when it hits $y$, and is independent of
$\left(R_{t};\;t\leq \rho\right)$; its infinitesimal generator is
given by:
$$\frac{1}{2}\dfrac{d^{2}}{dx^{2}}+\left(c\left(x\right)+\mathbf{1}_{\left(x>y\right)}\dfrac{s'\left(x\right)}{s\left(x\right)-s\left(m\right)}\right)\dfrac{d}{dx}.$$
\end{enumerate}
\end{itemize}
\end{prop}
\begin{proof}
The proof follows exactly the same lines as the proof of Proposition
\ref{decodifrec} and so we will not give it.
\end{proof}

\section*{Acknowledgements}
I am deeply indebted to Marc Yor whose ideas and help were essential
for the development of this paper. I would also like to thank Carl
Mueller for helpful conversations.

\end{document}